\documentclass[11pt]{article}
\usepackage{amsmath}
\usepackage{amsthm}
\usepackage{amsfonts}
\usepackage{amssymb}
\usepackage{latexsym}
\usepackage{fullpage}
\usepackage{mathrsfs}
\usepackage{mathcomp}
\usepackage{textcomp}
\usepackage{pifont}

\usepackage[dvips]{graphicx}

\newcounter{def_counter}

\newcounter{prob_counter}
\newcounter{as_counter}

\theoremstyle{plain}
\newtheorem{thm}{Theorem}
\newtheorem{cor}[thm]{Corollary}
\newtheorem{lem}[thm]{Lemma}
\newtheorem{prop}[thm]{Proposition}

\newtheorem{problem}[prob_counter]{Problem}

\theoremstyle{definition}
\newtheorem{defn}[def_counter]{Definition}

\newtheorem{assumption}[as_counter]{Assumption}

\newcommand{\N}{\mathbb N}

\newcommand{\mb}[1]{\mathbf{#1}}

\newcommand{\w}{\mathbf{w}}

\newcommand{\subword}{\prec}  
\newcommand{\rs}{\rho} 
\newcommand{\ls}{\sigma} 
\newcommand{\rp}{p_R} 
\newcommand{\lp}{p_L} 

\newcommand{\comp}[1]{\overline{#1}} 

\DeclareMathOperator{\alp}{alph}%
\DeclareMathOperator{\SW}{Sub} 
\DeclareMathOperator{\OC}{Occ} 

\title{Morphic and Automatic Words: Maximal Blocks and Diophantine Approximation}
\author{Yann Bugeaud\\
Universit\'e de Strasbourg\\  
U. F. R. de math\'ematiques\\  
7 rue Ren\'e Descartes, 67084 Strasbourg\\
France\\
{\tt bugeaud@math.u-strasbg.fr}\\ \\
Dalia Krieger\\
Faculty of Mathematics and Computer Science\\
The Weizmann Institute of Science\\
POB 26, Rehovot 76100\\
Israel\\
{\tt dalia.krieger@gmail.com}\\ \\
Jeffrey Shallit\\
School of Computer Science\\
University of Waterloo\\
Waterloo, ON N2L 3G1\\
Canada\\
{\tt shallit@cs.uwaterloo.ca}}

\begin {document}
\maketitle

\begin{abstract}
Let $\mb w$ be a morphic word over a finite alphabet $\Sigma$, and let $\Delta$ be a nonempty
subset of $\Sigma$. We study the behavior of maximal blocks consisting only of letters from
$\Delta$ in $\mb w$, and prove the following: let $(i_k,j_k)$ denote the starting and ending
positions, respectively, of the $k$'th maximal $\Delta$-block in $\mb w$. Then
$\limsup_{k\rightarrow\infty} (j_k/i_k)$ is algebraic if $\mb w$ is morphic, and rational if $\mb
w$ is automatic. As a result, we show that the same conclusion holds if $(i_k,j_k)$ are the
starting and ending positions of the $k$'th maximal zero block, and, more generally, of the $k$'th
maximal $x$-block, where $x$ is an arbitrary word.
This enables us to draw conclusions about the irrationality exponent of automatic and morphic
numbers. In particular, we show that the irrationality exponent of automatic (resp., morphic)
numbers belonging to a certain class that we define is rational (resp., algebraic).
\end{abstract}

\section{Introduction}\label{sect:intro}
The irrationality exponent $\mu(\xi)$ of an irrational number $\xi$ is the supremum of the real
numbers $\mu$ such that the inequality
$$
\biggl| \xi - {p \over q} \biggr| < {1 \over q^{\mu}}
$$
has infinitely many solutions in rational numbers $p/q$.

It follows from the theory of continued fractions that the irrationality exponent of every
irrational number is greater than or equal to $2$, and from the Borel-Cantelli lemma that it is
precisely equal to $2$ for almost all real numbers (with respect to Lebesgue measure). However, to
determine the irrationality exponent of a given real number $\xi$ is often a very difficult
problem. Na\"\i vely, we could hope to be able to read it off from the expansion of $\xi$ in some
integer base $b$, but this is almost never the case (see \cite{Amou&Bugeaud:2009} for a thorough
discussion).

Let $b \ge 2$ be an integer. Recently, Bugeaud \cite{Bugeaud:2008}
constructed a class
of real numbers whose irrationality exponent can be read off from their base-$b$ expansion.  This
class includes numbers of the form
$$
\xi_{\bf n} = \sum_{j \ge 1} \, b^{-n_j},
$$
for a strictly increasing sequence ${\bf n} = (n_j)_{j \ge 1}$ of positive
integers satisfying $n_{j+1} / n_j \ge 2$, for $j \ge 1$. To obtain good rational approximations to
$\xi_{\bf n}$, we simply truncate the above sum. Thus, we set
$$
\xi_{{\bf n}, J} = \sum_{j = 1}^J  \, b^{-n_j} = {p_J \over b^{n_J}}, \quad J \ge 1.
$$
It then follows from
$$
\biggl|\xi_{\bf n} - {p_J \over b^{n_J}} \biggr| < {2 \over (b^{n_J})^{n_{J+1}/n_J}}
$$
that
$$
\mu (\xi_{\bf n}) \ge \limsup_{j \to \infty} \, n_{j+1}/n_j.
$$
Shallit \cite{Shallit:1982} proved that the
continued fraction expansions of any
such $\xi_{\bf n}$ can be given explicitly,
and Bugeaud \cite{Bugeaud:2008}
proved its irrationality exponent is given by
\begin{equation}\label{eq:iratExp}
\mu (\xi_{\bf n}) = \limsup_{j \to \infty} \, n_{j+1}/n_j,
\end{equation}
and hence can be read off from its base-$b$ expansion. This simply means that the best rational
approximations to $\xi_{\bf n}$ are obtained by truncating its base-$b$ expansion.  In this
paper, among other results,
we use this method to study the irrationality exponents
of \emph{automatic and morphic numbers}. We let ${\cal C}$ denote
the class of all real numbers $\xi_{\bf n}$ as above.

Pure morphic words are infinite words generated by iterating a morphism defined over a finite
alphabet. Let $\Sigma$ be a finite alphabet, let $\epsilon$ denote the empty word, and let
$h:\Sigma^*\rightarrow\Sigma^*$ be a morphism. If there exists a letter $a\in\Sigma$ such that
$h(a) = ax$ for some $x\in\Sigma^+$, and furthermore, $h^n(a) \neq \epsilon$ for all $n\geq 0$ ($h$
is \emph{prolongable} on $a$), then the sequence $a, h(a), h^2(a),\ldots$ converges as $n$ tends to
infinity to the infinite word $h^\omega(a) = axh(x)h^2(x)\cdots$, which is a fixed point of the
morphism $h$. Such infinite fixed points are called \emph{pure morphic words}. An infinite word is
\emph{morphic} if it is the image under a \emph{coding} (that is, a letter-to-letter morphism) of a
pure morphic word; it is \emph{automatic} if it is morphic, and the underlying pure morphic word
can be generated by a \emph{uniform} morphism, that is, a morphism that maps all letters to words
of equal length. (Note: the standard definition of automatic words, or sequences, uses finite
automata. See, e.g., \cite[Chapter 5]{Allouche&Shallit:2003}.) A real number is \emph{automatic}
(resp., \emph{morphic}) if its expansion in some integer base $b \ge 2$ is an automatic (resp.,
morphic) word over the alphabet $\Sigma_b = \{0, 1, \ldots , b-1\}$.

Recall that a Liouville number is a real number $\xi$ satisfying $\mu(\xi) = \infty$. Adamczewski
and Cassaigne \cite{Adamczewski&Cassaigne:2006} established in 2006 that the irrationality exponent
of an automatic number is always finite, that is, automatic numbers are not Liouville numbers. In
\cite{Bugeaud:2008}, Bugeaud used Eq.~(\ref{eq:iratExp}) to show that any rational number $\mu \ge
2$ is the irrationality exponent of some automatic number. These two results motivate the following
question:

\begin{problem}\label{prob:automatic}
Determine the set of irrationality exponents of automatic numbers.
In particular, is the irrationality exponent of an automatic number always rational?
\end{problem}

Unfortunately, we are unable to settle this problem. However, our Theorem~\ref{thm:main} below
implies that the irrationality exponent of any automatic number in the class ${\cal C}$ is
rational. Consequently, to give a negative answer to Problem 1 we would need to use a radically
different method to construct automatic numbers.

As automatic numbers form a subclass of the morphic numbers, Problem~\ref{prob:automatic} can be
naturally extended as follows:

\begin{problem}\label{prob:morphic}
Determine the set of irrationality exponents of morphic numbers.
\end{problem}

Theorem~\ref{thm:main} implies that the
irrationality exponent of every morphic number in the class ${\cal C}$
is always an algebraic number.
Moreover, using the same method as in \cite{Bugeaud:2008}, we are able to show that every Perron
number $\mu \ge 2$ is the irrationality exponent of some morphic number (recall that a \emph{Perron
number} is a positive real algebraic integer that is greater in absolute value than all of its
conjugates). However, Problem~\ref{prob:morphic} remains unsolved.

Theorem~\ref{thm:main} is proved through a combinatorial study of maximal $\Delta$-blocks in
automatic and morphic words, where $\Delta\subseteq\Sigma$ is a subalphabet. We find the
combinatorial results interesting in their own right.

Our paper is organized as follows. In Section~\ref{sect:defn} we give some definitions and state
the main theorems, as well as some open problems. In Section~\ref{sect:pure} we analyze the
structure of $\Delta$-blocks in pure morphic words. In Section~\ref{sect:morphic} we apply the
results of Section~\ref{sect:pure} to morphic words in general. In Section~\ref{sect:perron} we
construct, for a given Perron number $\mu \geq2$, a morphic number $\xi\in\cal C$ such that
$\mu(\xi) = \mu$.

\section{Exponents of Diophantine approximation and maximal blocks}\label{sect:defn}

To carefully investigate the question whether one can read off the irrationality exponent of a real
number from its expansion in some integer base, Amou and Bugeaud \cite{Amou&Bugeaud:2009}
introduced new exponents of Diophantine approximation. Throughout the present paper, $|| \cdot ||$
denotes the distance to the nearest integer and $\lfloor  \cdot \rfloor$ denotes the greatest
integer function.

\begin{defn}
Let $\xi$ be an irrational real number. Let $b$ be an integer with $b \ge 2$. We let $v_b (\xi)$
denote the supremum of the real numbers $v$ for which the equation
$$
|| b^n \xi || < (b^n)^{-v}
$$
has infinitely many solutions in positive integers $n$. We let $v'_b (\xi)$ denote the supremum of
the real numbers $v$ for which the equation
$$
|| b^r (b^s - 1) \xi || < (b^{r+s})^{-v}
$$
has infinitely many solutions in positive integers $r$ and $s$.
\end{defn}

The exponent $v_b$ measures the accuracy with which a real number is approximable by rationals
obtained by truncating its base-$b$ expansion, while $v'_b$ measures the accuracy with which a real
number is approximable by rationals obtained by truncating its base-$b$ expansion and completing by
periodicity.

For every irrational number $\xi$, we have $v'_b (\xi) \ge v_b (\xi) \ge 0$ for $b \ge 2$, and
$$
\mu(\xi) \ge 1 + \max\{v'_b(\xi), 1\} \ge1 + \max\{v_b(\xi), 1\}, \quad \hbox{for $b \ge 2$}.
$$
Furthermore, any real number
$$
\xi_{\bf n} = \sum_{j \ge 1} \, b^{-n_j}
$$
belonging to the class ${\cal C}$ satisfies
\begin{equation} \label{eq:vb}
\mu (\xi) = 1 + v_b (\xi) = 1 + v'_b (\xi) =  \limsup_{j \to \infty} \, n_{j+1}/n_j.
\end{equation}
To understand the shift by $1$, just observe that $\mu (\xi) - 1$ is the supremum of the real
numbers $\mu$ such that
$$
|| q \xi || < q^{-\mu}
$$
has infinitely many solutions in positive integers $q$.

The main result of the present paper is the following:

\begin{thm}\label{thm:main}
Let $\xi$ be an irrational
real number, and suppose the expansion of $\xi$ in some integer base $b \geq 2$ is an automatic
(resp., morphic) word over
the alphabet $ \{0, 1, \ldots , b-1\}$. Then the
number $v_b (\xi)$ is finite and rational (resp., algebraic).
\end{thm}

Conversely, we do not know whether for every positive algebraic number $v$ there exist $b \ge 2$
and a morphic number $\xi$ such that $v_b (\xi) = v$. The next theorem provides a partial result
towards the resolution of this problem.

\begin{thm}\label{thm:perron}
For every rational number $v \ge 1$ (resp., Perron number $v > 1)$ and every integer $b \geq 2$
there exists a real number $\xi$, such that the base-$b$ expansion of $\xi$ is an automatic (resp.,
morphic) word over the alphabet $\{0, 1, \ldots , b-1\}$ and $v_b(\xi) = v - 1$.
\end{thm}

Theorem \ref{thm:perron} asserts that the set of values taken by the 
exponent $v_b$ at automatic irrational real numbers   
is precisely the set of nonnegative rational numbers.  

The real numbers $\xi$ constructed in the proof of Theorem \ref{thm:perron} satisfy $\mu (\xi) = 1
+ v_b (\xi)$ when $v_b(\xi) \geq 1$, but we do not know their irrationality exponent if $v_b (\xi)$
is less than $1$.

We use the same method as in \cite{Bugeaud:2008} to show the following:

\begin{cor}\label{cor:expmorphic}
For every Perron number $\mu \geq 2$ there exists a morphic number $\xi$ such that  $\mu (\xi) =
\mu$.
\end{cor}

We are currently unable to determine the set of positive algebraic numbers $v$ for which there
exist $b \ge 2$ and a morphic number $\xi$ such that $v_b (\xi) = v$. But this set strictly
contains the union of the positive rational numbers and the numbers of the form $r-1$ with $r$ a
Perron number.

Theorem \ref{thm:main} and Corollary \ref{cor:expmorphic} can be phrased in combinatorial terms.
Let $\w = w_0w_1w_2\cdots$ be an infinite word over $\Sigma_b = \{0, 1, \ldots , b-1\}$, and let $0
\leq i \leq j$. We say that $w_i\cdots w_j$ is a \emph{maximal zero block} in $\w$ if $w_i =
w_{i+1} = \cdots = w_j = 0$, $w_{j+1} \neq 0$, and either $i = 0$ or $w_{i-1} \neq 0$. Theorem
\ref{thm:main} and Corollary~\ref{cor:expmorphic} can be phrased in terms of the maximal zero
blocks in the base-$b$ expansion of $\xi$:

\begin{thm}\label{thm:maina}
Let $\mb w = w_0w_1w_2\cdots$ be an automatic (resp., morphic) word over $\{0, 1, \ldots , b-1\}$,
that does not have a suffix of the form $0^\omega$. For $k\geq 0$, let $(i_k, j_k)$ denote the
starting and ending positions, respectively, of the $k$'th maximal zero block in $\mb w$. Then
$\limsup_{k\rightarrow\infty} j_k/i_k$ is finite and rational (resp. algebraic).
\end{thm}

Strictly speaking, Theorem \ref{thm:main} is not a restatement of Theorem \ref{thm:maina}, since,
to deal with the exponent $v_b$, we also have to control the occurrences of blocks composed only of
the digit $b-1$.

\begin{thm}\label{thm:perrona}
For every Perron number $\mu \ge 2$ there exist a binary morphic word $\mb w = w_0w_1w_2\cdots$,
such that the sequence of indices $\{n_j\}_{j\geq 0} = \{n: w_n = 1\}$ satisfies
\begin{enumerate}
    \item $n_{j+1}/n_j \geq 2$ for all $j \geq 0$;
    \item $\limsup_{j\to\infty} n_{j+1}/n_j =  \mu$.
\end{enumerate}

\end{thm}

Maximal zero blocks are a special case of \emph{maximal $x$-blocks}, where $x\in\Sigma^+$ is an
arbitrary word. We say that $y = w_i\cdots w_j$ in $\w$ is an $x$-block if there exist some proper
suffix $x'$ and proper prefix $x''$ of $x$, such that $y = x'x^nx''$ for some integer $n \geq 1$;
if $x''w_{i+1}$ is not a prefix of $x$, and either $i = 0$ or $w_{i-1}x'$ is not a suffix of $x$,
then the $x$-block is maximal. For example, let $\mb w = 0100111010101000\cdots\in\{0,1\}^\omega$,
and let $x = 01$. Then $w_0w_1w_2 = 010$,  $w_3w_4 = 01$, and $w_6\cdots w_{13} = 10101010$ are all
maximal $x$-blocks. Theorem~\ref{thm:maina} can be generalized as follows:

\begin{thm}\label{thm:maina2}
Let $\mb w = w_0w_1w_2\cdots$ be an automatic (resp., morphic) word over $\{0, 1, \ldots , b-1\}$,
and let $x\in\Sigma_b^+$. Assume $\w$ does not have a suffix of the form $x^\omega$. For $k\geq 0$,
let $(i_k, j_k)$ denote the starting and ending positions, respectively, of the $k$'th maximal
$x$-block in $\mb w$. Then $\limsup_{k\rightarrow\infty} j_k/i_k$ is finite and rational (resp.,
algebraic).
\end{thm}

Theorems~\ref{thm:maina} and~\ref{thm:maina2} are proved in Section~\ref{sect:morphic}.
Theorems~\ref{thm:perrona} and~\ref{thm:perron} and Corollary~\ref{cor:expmorphic} are proved in
Section~\ref{sect:perron}.

In light of Theorem~\ref{thm:maina2}, it seems plausible that the method used in the proof would
allow us to say something about the exponent $v'_b$. This is not the case, however, since we then
have to consider possible repetitions of every word $x\in\Sigma_b^+$, and hence, to take the
supremum of an infinite set of rational (resp., algebraic) numbers, that may not be all distinct.
We cannot guarantee that this supremum is rational (resp., algebraic), nor even that it is finite.
Thus, we are unfortunately unable to establish the following
statement:\\

\noindent \textbf{Unproven Assertion (i)}.
\emph{Let $\xi$ be an irrational
real number, and suppose the expansion of $\xi$ in some integer base $b \geq 2$ is an automatic
(resp., morphic) word over the alphabet $\Sigma_b = \{0, 1, \ldots , b-1\}$. Then the number $v'_b
(\xi)$ is finite and rational
(resp., algebraic)}.\\

If we could prove that for every irrational
real number $\xi$, such that the expansion of $\xi$ in some integer base $b\geq 2$ is morphic
and has
sublinear complexity, the number $v'_b (\xi)$ is finite, then we could extend Theorem 2.1 of
Adamczewski and Cassaigne \cite{Adamczewski&Cassaigne:2006}, asserting that the irrationality
exponent of an automatic number is always finite,
as follows:\\

\noindent \textbf{Unproven Assertion (ii)}.
\emph{A morphic number of sublinear complexity cannot be
a Liouville number}.\\

Lemma 5.1 of \cite{Adamczewski&Cassaigne:2006} states that $v'_b (\xi)$ is finite for every
irrational automatic number, and is a key step in the proof of Theorem 2.1 of
\cite{Adamczewski&Cassaigne:2006}. Here, the assumption that $\xi$ is automatic is crucial. The
other steps of the proof do not require such a strong condition on $\xi$ and can be easily adapted
to the case where $\xi$ is morphic with sublinear complexity.

However, for numbers in the class $\cal C$, Theorem~\ref{thm:main} and Equality~(\ref{eq:vb}) imply
the following corollary:
\begin{cor}\label{cor:main}
Let $\xi$ be an automatic (resp., morphic) number in the class $\cal C$.
Then $\mu(\xi)$ is finite and rational (resp., algebraic).
\end{cor}

Corollary \ref{cor:main} is a small step towards the resolution of Problem~\ref{prob:automatic}.

\section{$\Delta$-blocks in pure morphic words}\label{sect:pure}

\begin{defn}\label{def:occ}
Let $\mb w = w_0w_1w_2\cdots\in\Sigma^\omega$. Let $\SW(\mb w)$ denote the set of finite subwords
of $\mb w$. An \emph{occurrence} of $\mb w$ is a triple $(u,i,j)$, where $\epsilon\neq u\in \SW(\mb
w)$ and $0\leq i\leq j$, such that $w_i\cdots w_j = u$. We usually denote an occurrence $(u,i,j)$
simply by $u$. The set of all occurrences of $\mb w$ is denoted by $\OC(\mb w)$. An occurrence
$(u,i,j)\in\OC(\mb w)$ \emph{contains} an occurrence $(u',i',j')\in\OC(\mb w)$, denoted $u'\subword
u$, if $i\leq i'$ and $j\geq j'$.
\end{defn}

\begin{defn}\label{def:delta-block}
Let $\mb w = w_0w_1w_2\cdots\in\Sigma^\omega$. Let $\Delta\subseteq\Sigma$, and let $\comp{\Delta}
= \Sigma\setminus\Delta$. An occurrence $(u,i,j)\in\OC(\mb w)$ is a \emph{$\Delta$-block} if
$u\in\Delta^+$. A $\Delta$-block $(u,i,j)\in\OC(\mb w)$ is \emph{maximal} if
$w_{j+1}\in\comp{\Delta}$, and either $i = 0$ or $w_{i-1}\in\comp{\Delta}$.
\end{defn}

Our goal in this section is to prove the following theorem:
\begin{thm}\label{thm:main_pure}
Let $h:\Sigma^*\rightarrow\Sigma^*$ be a nonerasing morphism, and let $\mb w = w_0w_1w_2\cdots =
h^\omega(w_0)$. Let $\Delta\subsetneq\Sigma$ be a nonempty subalphabet, such that $\mb w$ contain
infinitely many letters of $\comp{\Delta}$ and $\Delta$-blocks of unbounded length. For $k =
0,1,2,\ldots$, let $(i_k,j_k)$ denote the starting and ending positions, respectively, of the
$k$'th maximal $\Delta$-block in $\mb w$. Then $\limsup_{k\rightarrow\infty} j_k/i_k$ is an
algebraic number of degree at most $|\Sigma|$. If $h$ is also uniform, then
$\limsup_{k\rightarrow\infty} j_k/i_k$ is rational.
\end{thm}

We require $\mb w$ to contain $\Delta$-blocks of unbounded length because otherwise
$\limsup_{k\rightarrow\infty} j_k/i_k$ is trivially rational. This condition implies in particular
that $\mb w$ is \emph{aperiodic}, that is, it is not \emph{ultimately periodic}. Here an ultimately
periodic word is a word of the form $\mb w = xy^\omega$, where $x\in\Sigma^*$ and $y\in\Sigma^+$.

Proving Theorem~\ref{thm:main_pure} will enable us to prove the algebraicity (resp., rationality)
of the sequence of zero blocks in morphic (resp., automatic) words in general: if $\mb w =
\tau(h^\omega(a))$, where $\tau$ is a coding, then a maximal zero block in $\mb w$ is the image
under $\tau$ of a maximal $\Delta$-block in $h^\omega(a)$, where $\Delta = \tau^{-1}(0)$. The case
of maximal $x$-blocks will be proved by applying a morphic-preserving (resp.,
automaticity-preserving) transformation to $\mb w$.

The technique we use to prove Theorem~\ref{thm:main_pure} is very similar to the technique used to
prove the algebraicity of critical exponents in pure morphic words
\cite{Krieger:2007,Krieger:thesis}. The idea is as follows:
\begin{enumerate}
  \item The sequence of maximal $\Delta$-blocks can be partitioned into subsequences, where for
  each subsequence, every element is an image under $h$ of the previous element, up to a small
  change at the edges.
  \item There are only finitely many different such subsequences in $\mb w$. Since we are
  interested in $\limsup$, it is enough to consider only the first of each of the different
  subsequences.
  \item The $\limsup$ of a subsequence can be computed using the incidence matrix of $h$ (see
  Definition~\ref{def:incidence}). In particular, the $\limsup$ is a rational
  expression of the eigenvalues of the said matrix, which are algebraic numbers of degree at most $|\Sigma|$.
  \item When $h$ is uniform, the expression turns out to be rational.
\end{enumerate}

\begin{defn} \label{def:incidence}
Let $\Sigma = \Sigma_n = \{0,1,\ldots,n-1\}$, let $h:\Sigma_n^*\rightarrow\Sigma_n^*$, and let
$u\in\Sigma_n^*$. The \emph{Parikh vector} of $u$, denoted by $[u]$, is a vector of size $n$ that
counts how many times different letters occur in $u$: $[u] = (|u|_0,|u|_1, \ldots ,|u|_{n-1})^T$.
The \emph{incidence matrix} associated with $h$, denoted by $A(h)$, is an $n\times n$ matrix, whose
$j$th column is the Parikh vector of $h(j)$:
$$
A(h) = (a_{i,j})_{0\leq i,j<n} \;;\;\; a_{i,j} = |h(j)|_i\;.
$$
\end{defn}

\begin{prop} \label{prop:incidence}
Let $h:\Sigma^*\rightarrow\Sigma^*$, and let $A = A(h)$. Then:
\begin{enumerate}
  \item $[h(u)] = A[u]$ for all $u\in\Sigma^*$;
  \item $A(h^n) = A^n$ for all $n\in\N$.
\end{enumerate}
\end{prop}

See, e.g., \cite[Section 8.2]{Allouche&Shallit:2003}.\\

\noindent \textbf{Notation:} for a word $w$ (finite or not), $\alp(w)$ denotes the set of letters
occurring in $w$.

\begin{lem}\label{lem:reachable}
Let $h:\Sigma^*\rightarrow\Sigma^*$. Then there exists some power $g$ of $h$ such that for all
$a\in\Sigma$ and for all $n \geq 1$, $\alp(g^n(a)) = \alp(g(a))$.
\end{lem}
\begin{proof}
Let $A = A(h) = (a_{i,j})$, and denote $A^n = A(h^n) = (a_{i,j}^{(n)})$. Then for all letters
$a,b\in\Sigma$ and for all $n\geq 1$, $b\in\alp(h^n(a))$ if and only if $|h^n(a)|_b > 0$, that is,
if and only if $a_{b,a}^{(n)} > 0$. Since we care only about the zero pattern of $A^n$ and not
about the value of the non-zero entries, it is enough to consider $A$ as a boolean matrix. Let $B$
be a $|\Sigma|\times|\Sigma|$ boolean matrix, such that $b_{i,j} = 0$ if and only if $a_{i,j} = 0$.
Then it is enough to prove the following: there exists some power $B'$ of $B$, such that $B'^n =
B'$ for all $n\geq 1$.

Since there are only finitely many boolean matrices of a given size, there exist some integers
$t\geq 0$ and $c\geq1$ such that $B^t = B^{t+c}$, and so $B^{t+k} = B^{t+k+nc}$ for all
$k\in\{0,1,\ldots,c-1\}$ and for all $n\geq 0$. Choose a $k\in\{0,1,\ldots,c-1\}$ such that
$c|t+k$, and let $B' = B^{t+k}$. Then $t+k = mc$ for some integer $m$, and for all $n \geq 1$,
$$
B'^n = B^{n(t+k)} = B^{t+k+(n-1)mc} = B^{t+k} = B'.
$$
By setting $g = h^{t+k}$ we get the desired morphism.
\end{proof}

Let $\mb w = h^\omega(a)$ be a pure morphic word over $\Sigma$. Then $\mb w = (h^t)^\omega(a)$ for
all $t \geq 1$, and so we can replace $h$ by some convenient power. Therefore, by
Lemma~\ref{lem:reachable}, we can assume the following:

\begin{assumption}\label{as:R}
For all $a\in \Sigma$ and for all $n\geq 1$, $\alp(h^n(a)) = \alp(h(a))$.
\end{assumption}

In addition, for the rest of this section we assume that $h$ is nonerasing.

\begin{defn}\label{def:inverse_img}
Let $h:\Sigma^*\rightarrow\Sigma^*$, and let $\mb w = w_0w_1w_2\cdots = h^\omega(w_0)$. The
\emph{inverse image} under $h$ of an occurrence $u\in\OC(\mb w)$, denoted $h^{-1}(u)$, is the
shortest occurrence $v\in\OC(\mb w)$ such that $h(v)$ contains $u$.
\end{defn}

Note that for an occurrence $u$ (rather than a subword $u$), the inverse image is well defined.
Whenever we use the notation $h^{-1}(u)$ it should be understood that $u$ is an occurrence.

In the next two lemmas, we want to establish the following idea: if a pure morphic word $\mb w =
h^\omega(w_0)$ contains $\Delta$-blocks of unbounded length and infinitely many letters of
$\comp{\Delta}$, then sufficiently long $\Delta$-blocks are images under $h$ of other
$\Delta$-blocks, except perhaps for edges of a bounded length.

\begin{lem}\label{lem:inverse_img}
Let $\mb w = w_0w_1w_2\cdots = h^\omega(w_0)$ be an aperiodic pure morphic word over an alphabet
$\Sigma$, where $h$ is nonerasing and satisfies Assumption~\ref{as:R}, and let $M =
\max\{|h(a)|:a\in\Sigma\}$. Let $\Delta\subsetneq\Sigma$ be a nonempty subalphabet, such that $\mb
w$ contains infinitely many letters of $\comp{\Delta}$ and $\Delta$-blocks of unbounded length. Let
$u = w_r\cdots w_s\in\OC(\mb w)$ be a maximal $\Delta$-block, such that $|u| > M^2$ and $r > M$,
and let $h^{-1}(u) = w_i\cdots w_j$. Then
\begin{enumerate}
  \item $w_{i+M}\cdots w_{j-M}$ is a (not necessarily maximal) $\Delta$-block;
  \item $w_{i-M+1}\cdots w_{i+M-1}$ contains a letter of $\comp{\Delta}$;
  \item $w_{j-M+1}\cdots w_{j+M-1}$ contains a letter of $\comp{\Delta}$.
\end{enumerate}
\end{lem}
\begin{proof}
Suppose there is a letter $a$ occurring in $u$ such that $b := h^{-1}(a)\in\comp{\Delta}$. Let $c =
h^{-1}(b)$ (since $a$ occurs at a position $k > M \geq |h(w_0)|$, $h^{-1}(b)$ is well-defined).
Then by Assumption~\ref{as:R}, $h^2(c)$ contains $b$, and so $h^2(c) = h^2(h^{-2}(a))$ is not
contained in $u$. This implies that $b$ (which is contained in $h(c)$) occurs at a distance of at
most $M$ from the edges of $h^{-1}(u)$. Therefore, $w_{i+M}\cdots w_{j-M}$ is a $\Delta$-block.

Now consider $w_{i-M+1}\cdots w_{i+M-1}$. Since $u$ is maximal, $d := w_{r-1}\in\comp{\Delta}$. Let
$e = h^{-2}(d)$ (again, $h^{-2}(d)$ is well-defined, since $r-1\geq M$). Then $h(e)$ contains $d$
by Assumption~\ref{as:R}. But $h(e)$ is contained in $w_{i-M+1}\cdots w_{i+M-1}$, and so
$w_{i-M+1}\cdots w_{i+M-1}$ contains a letter of $\comp{\Delta}$. Similarly, $w_{j-M+1}\cdots
w_{j+M-1}$ contains a letter of $\comp{\Delta}$.
\end{proof}

\begin{lem}\label{lem:img}
Under the conditions of Lemma~\ref{lem:inverse_img}, let $u = w_i\cdots w_j\in\OC(\mb w)$ be a
maximal $\Delta$-block, such that $|u| > M^2$ and $i > M$. Then
\begin{enumerate}
  \item $h(w_{i+M}\cdots w_{j-M})$ is a $\Delta$-block;
  \item $h(w_{i-M+1}\cdots w_{i+M-1})$ contains a letter of $\comp{\Delta}$;
  \item $h(w_{j-M+1}\cdots w_{j+M-1})$ contains a letter of $\comp{\Delta}$.
\end{enumerate}
\end{lem}
\begin{proof}
Suppose $h(u)$ contains a letter $b\in\comp{\Delta}$. Then there exists a letter $a\in\Delta$ such
that $h(a)$ contains $b$. By the same argument as in the proof of Lemma~\ref{lem:inverse_img},
$h(h^{-1}(a))$ contains $b$, and cannot be contained in $u$. We get that $a$ occurs at a distance
of at most $M$ from the edges of $u$, and so $h(w_{i+M}\cdots w_{j-M})$ is a $\Delta$-block. The
rest is proved similarly.
\end{proof}

\begin{cor}\label{cor:sequences}
Under the conditions of Lemma~\ref{lem:inverse_img}, the set of maximal $\Delta$-blocks $u =
w_i\cdots w_j$ that satisfy $i > M$ and $|u|>M^2$ can be partitioned into (infinitely many)
sequences, each of which has the form $u^{(0)}, u^{(1)}, u^{(2)},\cdots$, where for all $k \geq 0$:
\begin{enumerate}
  \item $u^{(k)}$ is a maximal $\Delta$-block;
  \item if $u^{(k)} = w_i\cdots w_j$, then $h(w_{i+M}\cdots w_{j-M}) \subword u^{(k+1)} \subword h(w_{i-M+1}\cdots w_{j+M-1})$.
\end{enumerate}
\end{cor}

\begin{defn}\label{defn:orderedDelta}
Let $\mb w = w_0w_1w_2\cdots = h^\omega(w_0)$ be an aperiodic pure morphic word over an alphabet
$\Sigma$, let $M = \max\{|h(a)|:a\in\Sigma\}$, and let $\Delta\subsetneq\Sigma$. A
\emph{$\Delta$-sequence} in $\mb w$ is a sequence $u^{(k)} = w_{i_k}\cdots w_{j_k}$ of maximal
$\Delta$-blocks, where for all $k\geq 0$
\begin{itemize}
  \item $i_k > M$;
  \item $|u^{(k)}| > M^2$; and
  \item $h(w_{i_k+M}\cdots w_{j_k-M}) \subword u^{(k+1)} \subword h(w_{i_k-M+1}\cdots w_{j_k+M-1})$.
\end{itemize}
\end{defn}

\begin{defn}\label{defn:stretch}
Let $\mb w = w_0w_1w_2\cdots = h^\omega(w_0)$ and let $\{u^{(k)}\}_{k\geq0}$ be a
$\Delta$-sequence. For $k \geq 0$, let $h(w_{i_k}) = w_{r_{k+1}}\cdots w_{s_{k+1}}$ and $h(w_{j_k})
= w_{m_{k+1}}\cdots w_{n_{k+1}}$. Then
\begin{itemize}
  \item $u^{(k+1)}$ is \emph{growing} on the left if $i_{k+1} < r_{k+1}$;
  \item $u^{(k+1)}$ is \emph{shrinking} on the left if $i_{k+1} > r_{k+1}$;
  \item $u^{(k+1)}$ is \emph{stationary} on the left if $i_{k+1} = r_{k+1}$.
\end{itemize}
Similarly, $u^{(k+1)}$ is growing on the right if $j_{k+1} > n_{k+1}$, shrinking on the right if
$j_{k+1} < n_{k+1}$, and stationary on the right if $j_{k+1} = n_{k+1}$. The \emph{left stretch} of
$u^{(k+1)}$, denoted by $\ls^{(k+1)}$, is the word that occurs between the left edge of $u^{(k+1)}$
and the left edge of $h(u^{(k)})$. That is, if $u^{(k+1)}$ is shrinking on the left, then
$\ls^{(k+1)} := w_{r_{k+1}}\cdots w_{i_{k+1}-1}$ (in this case we say that the left stretch is
\emph{negative}); if $u^{(k+1)}$ is growing on the left, then $\ls^{(k+1)} := w_{i_{k+1}} \cdots
w_{r_{k+1}-1}$ (in this case we say that the left stretch is \emph{positive}). Note that if
$\ls^{(k+1)}$ is positive then it is contained in $u^{(k+1)}$, and if it is negative then it
borders $u^{(k+1)}$ on the left. If $r_{k+1} = i_{k+1}$ then $\ls^{(k+1)} := \epsilon$. The
\emph{right stretch}, denoted by $\rs^{(k)}$, is defined similarly.

The $k$'th \emph{left pivot}, denoted by $\lp^{(k)}$, is the rightmost letter in $w_{i_k-M+1}\cdots
w_{i_k+M-1}$ such that $h(\lp^{(k)})$ contains a letter of $\comp{\Delta}$; that is, $\lp^{(k)} :=
h^{-1}(w_{i_{k+1}-1})$. The \emph{right pivot}, denoted by $\rp^{(k)}$, is defined similarly.
\end{defn}

\begin{figure}[h]
\begin{center}
  \includegraphics [width=1\textwidth]{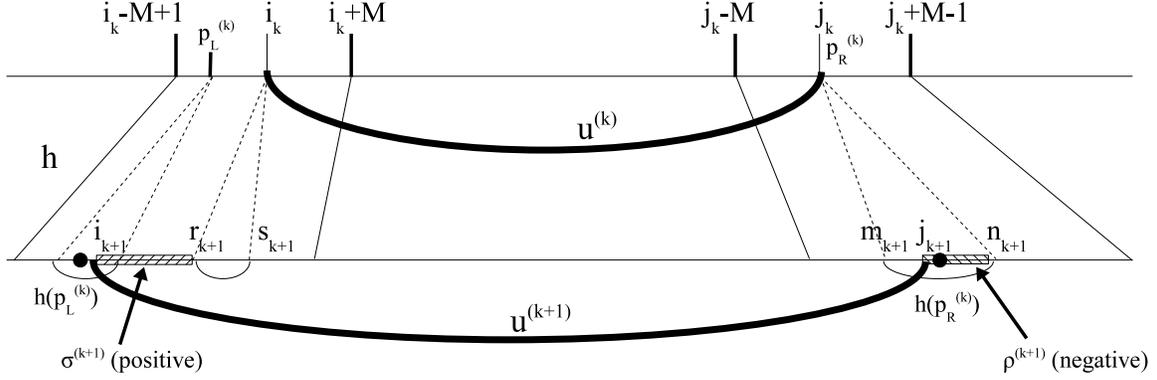}
  \caption{Maximal $\Delta$-blocks. $u^{(k+1)}$ is growing on the left and shrinking on the right.
  The $k$'th right pivot is exactly $w_{j_k}$. The black circles are the $\comp{\Delta}$ letters that terminate
  $u^{(k+1)}$.}
  \label{fig:deltablock}
\end{center}
\end{figure}

Figure~\ref{fig:deltablock} illustrates Definition~\ref{defn:stretch}.

\begin{lem}\label{lem:periodic}
Let $\mb w = w_0w_1w_2\cdots = h^\omega(w_0)$, where $h$ is nonerasing and satisfies
Assumption~\ref{as:R}, and let $\{u^{(k)}\}_{k\geq0}$ be a $\Delta$-sequence. Then the sequences
$\{\ls^{(k)}\}_{k\geq 1}$, $\{\rs^{(k)}\}_{k\geq 1}$ are ultimately periodic.
\end{lem}
\begin{proof}
We prove the lemma for the left stretch. The proof for the right stretch is similar.

For $k\geq 0$, let $\lp^{(k)} = w_{p_k}$. Then $h(w_{p_0+1}\cdots w_{i_0+M-1}) \in\Delta^+$, and so
$h^n(w_{p_0+1}\cdots w_{i_0+M-1}) \in\Delta^+$ for all $n > 0$ (recall Assumption~\ref{as:R}). In
particular, $h^2(w_{p_0+1}\cdots w_{i_0+M-1}) \in\Delta^+$, and so $\lp^{(1)}$ cannot occur in
$h(w_{p_0+1}\cdots w_{i_0+M-1})$. On the other hand, $\alp(h^2(w_{p_0})) = \alp(h(w_{p_0}))$, and
so $h(w_{p_0})$ contains a letter $a$ such that $h(a)$ contains a letter of $\comp{\Delta}$. In
particular, $\lp^{(1)} \subword h(\lp^{(0)})$. More generally, for all $k > 0$, $\lp^{(k)}$ is the
rightmost letter $a$ in $h(\lp^{(k-1)})$ such that $h(a)$ contains a letter of $\comp{\Delta}$.
This implies that the sequence of left pivots, $\{\lp^{(k)}\}_{k\geq 0}$, is ultimately periodic:
since $\Sigma$ is finite, there exist some $k\neq m$ such that $\lp^{(k)} = \lp^{(m)}$, and so
$\lp^{(k+n)} = \lp^{(m+n)}$ for all $n\geq 1$.

Now consider the left stretch. Let $h(\lp^{(0)}) = w_{x_1}\cdots w_{y_1}$. By definition, either
$i_1 = y_1+1$ (if the rightmost $\comp{\Delta}$ letter is the last letter of $h(\lp^{(0)})$), or
$x_1 < i_1 \leq y_1$. In the first case, $w_{i_1}\subword h(w_{p_0+1})$, and so
$h(w_{i_1})\in\Delta^+$, and $u^{(2)}$ is either growing or stationary. Since $h(\lp^{(0)})$
contains $\lp^{(1)}$ (that is, $h^2(\lp^{(0)})$ contains a letter of $\comp{\Delta}$), we get that
$\ls^{(2)}$ is a suffix of $h^2(\lp^{(0)})$.

If $x_1 < i_1 \leq y_1$, then $u^{(2)}$ can be also shrinking. However, in this case both
$h(w_{i_1})$ and $w_{i_2}$ are contained in $h^2(\lp^{(0)})$, and so $\ls^{(2)}\subword
h^2(\lp^{(0)})$. Similarly, $\ls^{(k)} \subword h(\lp^{(k-2)})$ for all $k \geq 2$. This implies
that $\{\ls^{(k)}\}_{k\geq 1}$ is ultimately periodic.
\end{proof}

Let $\{u^{(k)}\}_{k\geq0}$ be a $\Delta$-sequence. Since both $\{\ls^{(k)}\}_{k\geq 1}$ and
$\{\rs^{(k)}\}_{k\geq 1}$ are ultimately periodic sequences, the sequence
$\{(\ls^{(k)},\rs^{(k)})\}_{k\geq 1}$ is also ultimately periodic. By ignoring the first few
elements of the sequence we can assume it is purely periodic; by replacing $h$ by $h^p$, where $p$
is the period, we can partition $\{(\ls^{(k)},\rs^{(k)})\}_{k\geq 1}$ into $p$ subsequences, where
each subsequence has period 1. We now compute $i_k$ and $j_k$ for a sequence of maximal
$\Delta$-blocks, assuming that $\ls^{(k)}$ and  $\rs^{(k)}$ are fixed.

\begin{lem}\label{lem:algebraic}
Let $\{u^{(k)}\}_{k\geq0}$ be a $\Delta$-sequence, and assume that $\ls^{(k)} = \ls$ and $\rs^{(k)}
= \rs$ for all $k \geq 0$. Let $A$ be the incidence matrix of $h$, and let $\mb 1$ be the all ones
vector of size $1\times|\Sigma|$. Then there exist integral vectors $U,V,X,Y$ of size
$|\Sigma|\times 1$, where $U$ and $V$ are nonnegative and nonzero, and a constant $c$, such that
for all $k \geq 0$,
\begin{itemize}
  \item $i_k = \mb 1 \left(A^kV + (\sum_{n=0}^{k-1}A^k)X\right)$;
  \item $j_k = i_k + \mb 1 \left(A^kU + (\sum_{n=0}^{k-1}A^k)Y\right) - 1$;
  \item $j_k/i_k < c$.
\end{itemize}
\end{lem}
\begin{proof}
Let $v^{(k)} = w_0\cdots w_{i_k-1}$. Then $i_k = |v^{(k)}| = \mb 1[v^{(k)}]$ and $j_k = i_k +
|u^{(k)}|-1 = i_k + \mb 1[u^{(k)}] - 1$. To compute $i_k$ and $j_k$ we need to compute $[u^{(k)}]$
and $[v^{(k)}]$.

Let $U = [u^{(0)}]$ and $V = [v^{(0)}]$. Since $i_0 > M$, $v^{(0)}$ is a nonempty word, and so both
$U$ and $V$ are nonnegative, nonzero vectors. Let us assume for the moment that $\rs = \epsilon$.
Depending on whether $\ls$ is positive or negative, there are two possible situations:
\begin{enumerate}
  \item $h(u^{(k)}) = \ls u^{(k+1)}$ ($\ls$ is negative);
  \item $\ls h(u^{(k)}) = u^{(k+1)}$ ($\ls$ is positive).
\end{enumerate}
Suppose $\ls$ is negative. Then $h(u^{(0)}) = \ls u^{(1)}$, $h^2(u^{(0)}) = h(\ls u^{(1)}) =
h(\ls)h(u^{(1)}) = h(\ls)\ls u^{(2)}$, and by induction, $h^k(u^{(0)}) =
h^{k-1}(\ls)h^{k-2}(\ls)\cdots h(\ls)\ls u^{(k)}$. By Proposition~\ref{prop:incidence}, we get that
for all $k\geq 0$,
$$
[u^{(k)}] = [h^k(u^{(0)})] - [h^{k-1}(\ls)] - [h^{k-2}(\ls)] - \cdots - [\ls] = A^kU - \left(
\sum_{\ell=0}^{k-1} A^\ell \right)[\ls].
$$
Now suppose that $\ls$ is positive. Then $u^{(1)} = \ls h(u^{(0)})$, $u^{(2)} = \ls h(u^{(1)}) =
\ls h(\ls) h^2(u^{(0)})$, and by induction, $u^{(k)} = \ls h(\ls)\cdots h^{k-1}(\ls) h^k(u^{(0)})$.
By Proposition~\ref{prop:incidence}, we get that for all $k\geq 0$,
$$
[u^{(k)}] = [h^k(u^{(0)})] + [h^{k-1}(\ls)] + [h^{k-2}(\ls)] + \cdots + [\ls] = A^kU +
\left(\sum_{\ell=0}^{k-1} A^\ell\right)[\ls].
$$

If $\rs \neq \epsilon$, then, depending on its sign, we get that $[u^{(k)}] = A^kU \pm
(\sum_{\ell=0}^{k-1} A^\ell)[\ls] \pm (\sum_{\ell=0}^{k-1} A^\ell)[\rs]$. Let $Y =
\pm[\ls]\pm[\rs]$. Then for all $k\geq 0$,
$$
[u^{(k)}] = A^kU + \left(\sum_{\ell=0}^{k-1} A^\ell\right)Y.
$$

Similarly, $[v^{(k)}] = A^kV + (\sum_{\ell=0}^{k-1} A^\ell)[\ls]$ if $\ls$ is negative, and
$[v^{(k)}] = A^kV - (\sum_{\ell=0}^{k-1} A^\ell)[\ls]$ if $\ls$ is positive (here the roles are
inverted: if $\ls$ is negative then it is positive with respect to $v^{(k)}$, and vice versa). Let
$X = \mp [\ls]$. Then for all $k\geq 0$,
$$
[v^{(k)}] = A^kV + \left(\sum_{\ell=0}^{k-1} A^\ell\right)X.
$$

It remains to show that the sequence $j_k/i_k$ is bounded by a constant. Let $w^{(k)} =
w_{j_k+1}\cdots w_{i_{k+1}-1}$. Since maximal $\Delta$-blocks are disjoint and separated by at
least one letter from $\comp{\Delta}$, $|w^{(k)}| \geq 1$ for all $k\geq 0$. Now, for all $k \geq
1$,
$$
\begin{array}{lllll}
|u^{(k)}| &=& |h(u^{(k-1)})| \pm |\ls| \pm |\rs|  & < & M|u^{(k-1)}| + 2M\;,\\
|v^{(k)}| &=& |v^{(k-1)}u^{(k-1)}w^{(k-1)}| & > & |u^{(k-1)}|\;,
\end{array}
$$
and so
$$
\frac{j_k}{i_k} = 1 + \frac{|u^{(k)}| - 1}{|v^{(k)}|} < 1 + M + \frac{2M}{|u^{(0)}|}\;.
$$
This completes the proof of the lemma.
\end{proof}

The following theorem was proved in \cite{Krieger:2007,Krieger:thesis}:
\begin{thm}\label{thm:matpower}
Let $A$ be an $n\times n$ nonnegative integral matrix with no zero columns, and let $U,V,W$ be
nonnegative integral column vectors of size $n$, with $W\neq 0$. Let
$$
\mathcal F(k) = \frac{\mb 1\left(A^kU + (\sum_{i=0}^{k-1}A^i)V\right)}{\mb
1\left(A^kW\right)}\;,\;\;k\geq 0\;.
$$
Then
\begin{enumerate}
    \item $\{\mathcal F(k)\}_{k\geq0}$ has finitely many accumulation points;
    \item if $\alpha$ is a finite accumulation point of $\mathcal{F}$, then $\alpha$ is a rational
    expression of the eigenvalues of $A$. In particular, $\alpha$ is algebraic of degree at most $n$.
\end{enumerate}
\end{thm}

The proof of Theorem~\ref{thm:matpower} can be adapted, with slight changes, to the case of the
sequence $\{j_k/i_k\}_{k\geq 0}$. Here we have a sequence of the form
$$
\frac{j_k}{i_k} = 1+\frac{|u^{(k)}|-1}{|v^{(k)}|} = 1+\frac{\mb 1\left(A^k U +
(\sum_{\ell=0}^{k-1}A^\ell)Y\right) - 1}{\mb 1\left(A^kV+
(\sum_{\ell=0}^{m-1}A^\ell)X\right)}\;,\;\;k\geq 0\;,
$$
where $A$ is a $|\Sigma|\times |\Sigma|$ nonnegative integral matrix with no zero columns (recall
that $h$ is nonerasing), $U$ and $V$ are nonnegative integral vectors, both nonzero, and $X$ and
$Y$ are integral vectors, with possibly negative entries. However, since both $|u^{(k)}|$ and
$|v^{(k)}|$ are tending to infinity as $k$ tends to infinity, both nominator and denominator are
always positive. In particular, the fact that $X$ and $Y$ may contain negative entries does not
alter the result. Also, Lemma~\ref{lem:algebraic} implies that all accumulation points are finite.

\begin{cor}\label{cor:algebraic}
Under the conditions of Lemma~\ref{lem:algebraic}, $\limsup_{k\rightarrow\infty} j_k/i_k$ is an
algebraic number of degree at most $|\Sigma|$.
\end{cor}

\begin{lem}\label{lem:rational}
Under the conditions of Lemma~\ref{lem:algebraic}, if $h$ is uniform then
$\lim_{k\rightarrow\infty} j_k/i_k$ exists and is rational.
\end{lem}
\begin{proof}
If $h$ is an $m$-uniform morphism, then $|h(w)| = m|w|$ for all $w\in\Sigma^*$. Let $u =
|u^{(0)}|$, $v = |v^{(0)}|$, $y = \pm|\ls|\pm|\rs|$, and $x = \mp|\ls|$. Then the expressions for
$|u^{(k)}|$ and $|v^{(k)}|$ are reduced to
\begin{eqnarray}
\nonumber |u^{(k)}| &=& m^k|u| + (\sum_{\ell=0}^{k-1} m^\ell)y \;=\; m^k|u| + y\frac{m^k-1}{m-1}\;,\\
\nonumber |v^{(k)}| &=& m^k|v| + (\sum_{\ell=0}^{k-1} m^\ell)x \;=\; m^k|v| + x\frac{m^k-1}{m-1}\;.
\end{eqnarray}
Therefore,
$$
\frac{j_k}{i_k} = 1+\frac{|u^{(k)}|-1}{|v^{(k)}|} =
1+\frac{m^k|u|+y\frac{m^k-1}{m-1}-1}{m^k|v|+x\frac{m^k-1}{m-1}}
\;\;\;\overrightarrow{k\rightarrow\infty}\;\;\; 1+\frac{(m-1)|u|+y}{(m-1)|v|+x}\;.
$$
\end{proof}

\begin{proof}[\textbf{Proof of Theorem~\ref{thm:main_pure}}]
Let $M = \max\{|h(a)|:a\in\Sigma\}$. Since we are interested in $\limsup j_k/i_k$, it is enough to
consider only $\Delta$-blocks of size larger than $M^2$ that occur at an index $i > M$. By
Corollary~\ref{cor:sequences}, these $\Delta$-blocks can be partitioned into sequences, where for
each sequence, an element is the image under $h$ of the previous element, save perhaps for edges of
a bounded length. Let $u = w_i\cdots w_j$ be the first element of such a sequence. Then $|u| > M^2$
and $i > M$, and so by Lemma~\ref{lem:inverse_img}, $v:=h^{-1}(u)$ is a maximal $\Delta$-block (up
to the edges); however, $|v| \leq M^2$, or it would be part of the sequence itself. Taking into
account the occurrences of size $M$ on both sides of $v$, we get that each sequence is uniquely
determined by a subword of $\mb w$ of length at most $M^2+2M$. Since there are only finitely many
such subwords, there are only finitely many different such sequences. To compute the $\limsup$, it
is enough to consider only the first of each of the different sequences, where $i_k$ is the
smallest. Therefore, we need to consider only finitely many sequences. Each sequence can be further
partitioned into finitely many subsequences, where for each of those, $\limsup j_k/i_k$ is
algebraic of degree at most $|\Sigma|$ (Lemma~\ref{lem:algebraic}). For uniform morphisms, $\limsup
j_k/i_k$ is rational (Lemma~\ref{lem:rational}).
\end{proof}

\section{$\Delta$-blocks and $x$-blocks in morphic words}\label{sect:morphic}

In this section we extend Theorem~\ref{thm:main_pure} to morphic words in general, as described in
the beginning of Section~\ref{sect:pure}. First, the next theorem shows that we lose no generality
by restricting ourself to nonerasing morphisms:
\begin{thm}[{\cite[Theorem 7.5.1]{Allouche&Shallit:2003}}]\label{thm:7.5.1}
Every pure morphic word is the the image under a coding of a pure morphic word generated by a
nonerasing morphism.
\end{thm}

\begin{thm}\label{thm:main_morphic}
Let $\mb w$ be a morphic word over a finite alphabet $\Sigma$. Let $\Delta\subset\Sigma$ be a
nonempty proper subalphabet, such that $\mb w$ contain infinitely many letters of $\comp{\Delta}$
and $\Delta$-blocks of unbounded length. For $k = 0,1,2,\ldots$, let $u^{(k)} = w_{i_k}\cdots
w_{j_k}$ be the $k$'th maximal $\Delta$-block in $\mb w$. Then $\limsup_{k\rightarrow\infty}
j_k/i_k$ is algebraic. If $\mb w$ is also automatic then  $\limsup_{k\rightarrow\infty} j_k/i_k$ is
rational.
\end{thm}
\begin{proof}
Since $\mb w$ is morphic, there exists some alphabet $\Sigma'$, a morphism
$h:\Sigma'^*\rightarrow\Sigma'^*$, and a coding $\tau:\Sigma'^*\rightarrow\Sigma^*$, such that $\mb
w = \tau(h^\omega(a))$ for some $a\in\Sigma'$. By Theorem~\ref{thm:7.5.1}, we can assume that $h$
is nonerasing. Let $\Delta' = \tau^{-1}(\Delta)$. Then every maximal $\Delta$-block in $\mb w$ is
the image under $\tau$ of a maximal $\Delta'$-block in $h^\omega(a)$, and every maximal
$\Delta'$-block in $h^\omega(a)$ is mapped by $\tau$ to a maximal $\Delta$-block in $\mb w$. The
result follows from Theorem~\ref{thm:main_pure}.
\end{proof}


\begin{proof}[\textbf{Proof of Theorem~\ref{thm:maina}}]
Set $\Delta = \{0\}$ and apply Theorem~\ref{thm:main_morphic}.
\end{proof}


\begin{proof}[\textbf{Proof of Theorem~\ref{thm:maina2}}]
Let $|x| = d$, and let $X = \{(u^{(k)},i_k,j_k)\in\OC(\mb w): k\geq 0\}$. First, we partition $X$
into $d$ subsequences, $X_0,\ldots X_{m-1}$, where
$$
X_m = \{(u,r,s)\in X: u = x'x^nx'',\; x' \textrm{ is a proper prefix of } x, \textrm{ and }
r+|x'|\equiv m \pmod d\}.
$$
That is, $X_i$ is the sequence of maximal $x$-block for which $x$ itself begins at an index
equivalent to $m\pmod d$. For $m = 0,1,\ldots,d-1$ and for $k = 0,1,2,\ldots$, let
$(i_{k,m},j_{k,m})$ denote the starting and ending positions of the $k$'th element of $X_m$. Then
\begin{equation}\label{eq:limsup}
\limsup_{k\rightarrow\infty}\frac{j_k}{i_k} =
\max\left\{\limsup_{k\rightarrow\infty}\frac{j_{k,m}}{i_{k,m}}: 0\leq m < d\right\}\;.
\end{equation}

Let $\Sigma' = \{[a_0\cdots a_{d-1}]:a_0\cdots a_{d-1}\in\SW(\mb w)\}$, and define $d$ infinite
words over $\Sigma'$ by
\begin{eqnarray}
\nonumber \mb w_0 &=& [w_0\cdots w_{d-1}][w_d\cdots w_{2d-1}]\cdots,\\
\nonumber \mb w_1 &=& [w_1\cdots w_d][w_{d+1}\cdots w_{2d}]\cdots,\\
\nonumber \vdots && \vdots \\
\nonumber \mb w_{d-1} &=& [w_{d-1}\cdots w_{2d-2}][w_{2d-1}\cdots w_{3d-2}]\cdots.
\end{eqnarray}
By \cite[Theorem 7.9.1]{Allouche&Shallit:2003}, if $\mb w$ is morphic then $\mb w_m$ is morphic for
all $m$; by \cite{Cobham:1972}, if $\mb w$ is automatic then $\mb w_m$ is automatic for all $m$.
Let $\Gamma = \Sigma\cup\{\alpha\}$, where $\alpha\notin\Sigma$, and define a $d$-uniform morphism
$\tau:\Sigma'^*\rightarrow\Gamma^*$ by
$$
\tau([a_0\cdots a_{d-1}]) = \left\{
\begin{array}{ll}
\alpha^d, &\textrm{ if } a_0\cdots a_{d-1} = x;\\
a_0\cdots a_{d-1}, &\textrm{ if } a_0\cdots a_{d-1} \neq x.
\end{array}
\right.
$$
Let $\mb v_m = \tau(\mb w_m)$, $m = 0,\ldots, d-1$. By \cite[Corollary 7.7.5, Corollary
6.8.3]{Allouche&Shallit:2003}, if $\mb w_m$ is morphic (resp., automatic), then so is $\mb v_m$.
Let $\Delta = \{\alpha\}$, and let $(r_{k,m},s_{k,m})$ denote the starting and ending positions of
the $k$'th maximal $\Delta$-block in $\mb v_m$. Then for all $k\geq 0$, $|r_{k,m} - i_{k,m}| < d$
and $|s_{k,m} - j_{k,m}| < d$, and so $\limsup_{k\rightarrow\infty} s_{k,m}/r_{k,m} =
\limsup_{k\rightarrow\infty} j_{k,m}/i_{k,m}$. By Theorem~\ref{thm:main_morphic},
$\limsup_{k\rightarrow\infty} s_{k,m}/r_{k,m}$ is algebraic (resp., rational) if $\mb v_m$ is
morphic (resp., automatic), and so $\limsup_{k\rightarrow\infty} j_{k,m}/i_{k,m}$ is algebraic
(resp., rational) if $\mb w$ is morphic (resp., automatic). By (\ref{eq:limsup}), the result
follows.
\end{proof}


\section{Perron numbers as irrationality exponents of morphic numbers}~\label{sect:perron}

\begin{proof}[\textbf{Proof of Theorem~\ref{thm:perrona}}]
Let $\mu > 1$ be a Perron number. Then there exists a primitive integral square matrix $A$, of size
$k\times k$ for some positive integer $k$, such that $r(A) = \mu$, where $r(A)$ is the
Perron-Frobenius eigenvalue of $A$ \cite[Theorem 11.1.4]{Lind&Marcus:1995}. We may assume that
$k\geq 2$: if $\mu$ is not integral then necessarily $k \geq 2$, and if $\mu$ is integral we can
set $A$ to be the $2\times 2$ matrix ${\mu-1 \;\; \mu-1 \choose 1\;\;\;\;\;\;1}$. Let $\Sigma =
\Sigma_k = \{0,1,\ldots,k-1\}$, and let $h:\Sigma^*\to\Sigma^*$ be a morphism such that $A(h) = A$.
Then
$$
\frac{|h^n(0)|}{|h^{n-1}(0)|} = \frac{\mb 1A^n[0]}{\mb 1A^{n-1}[0]}\;,\;\;n\geq 0\;.
$$
Since $\mu$ is a Perron number, the Jordan decomposition of $A$ has one block of size 1 associated
with $\mu$, and it is easy to check that
$$
\lim_{n\to\infty}\frac{|h^n(0)|}{|h^{n-1}(0)|} = \mu.
$$
For $\mu>2$, the above equation implies that $|h^n(0)|/|h^{n-1}(0)| > 2$ for $n$ sufficiently
large; for $\mu=2$, we let $A = {1\;\;1\choose1\;\;1}$, and get that $|h^n(0)|/|h^{n-1}(0)| = 2$
for all $n$. Note that we do not require $h$ to be prolongable on $0$; we consider only the finite
words $\{h^n(0)\}_{n\geq0}$.

Let $\Gamma = \{\alpha,\beta\}\cup\Sigma$. Define a morphism $g:\Gamma^*\to\Gamma^*$ by $g(\alpha)
= \alpha\beta0$, $g(\beta) = \beta$, and $g(i) = h(i)$ for $i\in\Sigma$. Let $\mb u =
g^\omega(\alpha)$. It is an easy induction to show that
$$\mb u = \alpha\cdot\beta\cdot0\cdot\beta\cdot h(0) \cdot \beta \cdot h^2(0) \cdot \beta \cdots.$$
Now let $\tau:\Gamma^*\to\{0,1\}^*$ be the coding that maps $\beta$ to $1$ and all other letters to
$0$, and let $\w = \tau(\mb u)$. Then
$$\w = 0\cdot1\cdot0^{x_0}\cdot1\cdot 0^{x_1} \cdot 1 \cdot 0^{x_2} \cdot 1 \cdots,$$
where $x_n = |h^n(0)|$. Thus, $\w$ is a morphic sequence that satisfies the conditions of
Theorem~\ref{thm:perrona}. This completes the proof of the theorem.
\end{proof}

\begin{proof}[\textbf{Proof of Theorem~\ref{thm:perron}
and of Corollary \ref{cor:expmorphic}}] For morphic numbers, let $\mu > 1$ be a Perron number, and
let $b\geq 2$ be an integer. Let $\mb w$ be as in the proof of Theorem~\ref{thm:perrona}, and let
$\{n_j\}_{j\geq 0} = \{n: w_n = 1\}$. We associate with $\mb w$ the real number
$$\xi_{\w} = \sum_{j\ge 0}  \, b^{-n_j}.$$
Since $\limsup_{j \to \infty} \, n_{j+1}/n_j = \mu$, we get that
$$v_b (\xi_{\w}) = \mu-1.$$
Also, if $\mu \geq 2$ then $\xi_{\w}$ belongs to class ${\cal C}$, and so $\mu(\xi_{\w}) = \mu$.
This follows from the fact that $n_{j+1} \geq 2 n_j$ for $j$ sufficiently large, and this enables
us to use the Folding Lemma as in \cite{Shallit:1982,Bugeaud:2008} to construct the continued
fraction expansion of a rational translate of $\xi_{\w}$.

Now consider automatic numbers. For $v = 1$, we can choose any morphic binary word that does not
contain unbounded $0$-blocks or $1$-blocks (e.g., the Thue--Morse word) to be the base-$b$
expansion of $\xi$, where $b\geq 2$ is any integer, and get that $v_b (\xi) = 0 = v-1$. Suppose $v
= p/q > 1$. We define a binary infinite word $\mb u = u_0u_1u_2\cdots$ by letting $u_n = 1$ if and
only if $n$ belongs to the set
$$
\bigcup_{h \geq 0} \{(p)p^h, (p+1)p^h, \ldots, (qp)p^h \}.
$$
Then $\mb u$ is $p$-automatic, because its $p$-kernel contains only two sequences, namely the
sequence $0^\omega$ and the sequence $0^{p} 1^{qp - p + 1} 0^\omega$ (see \cite[Theorem
6.6.2]{Allouche&Shallit:2003}). Let $\{n_j\}_{j\geq 0} = \{n: u_n = 1\}$. Since $(p + 1)/p < p/q$,
we have
$$
\limsup_{j\to\infty} \frac{n_{j+1}}{n_j} = \lim_{h\to\infty} \frac{p\cdot p^{h+1}}{qp\cdot p^{h}}
=\frac pq\;.
$$
Consequently, the real number
$$\xi_{\mb u} = \sum_{j \ge 1} \, b^{-n_j}$$
satisfies
$$v_b (\xi_{\mb u}) = \frac pq - 1.$$
\end{proof}

\noindent \textbf{\emph{Remark 1}}. We stress that, with the above construction for automatic
numbers, we do not know the value of $\mu (\xi_{\mb u})$, because the condition $n_{j+1}/n_j \geq
2$ for all $j$ sufficiently large is not satisfied. This is not the case with the slightly more
complicated construction given in \cite{Bugeaud:2008}, which works under the assumption that $p/q$
exceeds $2$.\\

\noindent \textbf{\emph{Remark 2}}. For morphic numbers, we can go a bit further. Keep the notation
of the proof of Theorem~\ref{thm:perrona}. Let $\Gamma = \{\alpha,\beta,\gamma\}\cup\Sigma$. Let
$a$ and $b$ be non-negative integers. Define a morphism $g:\Gamma^*\to\Gamma^*$ by $g(\alpha) =
\alpha \beta0 (\gamma 0)^a \beta 0  (\gamma 0)^b$, $g(\beta) = \beta$, $g(\gamma) = \gamma$ and
$g(i) = h(i)$ for $i\in\Sigma$. Let $\mb u = g^\omega(\alpha)$. It is an easy induction to show
that
$$
\mb u = \alpha\cdot\beta\cdot0\cdot (\gamma \cdot 0)^a \cdot \beta \cdot 0 \cdot (\gamma \cdot 0)^b
\cdot \beta\cdot h(0) \cdot (\gamma \cdot h(0) )^a  \cdot \beta \cdot h(0)  \cdot (\gamma \cdot
h(0) )^b  \cdot \beta  \cdots.
$$
Now let $\tau:\Gamma^*\to\{0,1\}^*$ be the coding that maps $\beta$ to $1$ and all other letters to
$0$, and let $\w = \tau(\mb u)$. With suitable choices for $a$ and $b$, the value of the exponent
$v_b$ at the corresponding morphic number is neither rational, nor a Perron number minus $1$.
However, we do not know whether every positive algebraic number can be attained.

\end {document}